\documentclass[8pt]{article}

\usepackage{amsmath}
\usepackage{amsthm}
\usepackage{amsxtra}
\usepackage{amsfonts,amssymb}
\newtheorem{Th}{Theorem}
\newtheorem{Prop}{Proposition}

\theoremstyle{definition}
\newtheorem{Def}[Prop]{Definition}

\newtheorem{Rem}{Remark}

\date{}
\title{Stable states and representations of the infinite symmetric group}

\author{A.~M.~Vershik\footnote{Partially supported by the RFBR grant
10-01-90411-Ukr-a.} 
\and
N.~I.~Nessonov\footnote{Partially supported by the National Academy of
Sciences of
Ukraine grant 2010-2011-284.}}

\begin{document}
\maketitle
\begin{flushright}UDC 517.986\end{flushright}
\section{Introduction}

The problem studied in this paper can be briefly formulated as
follows: what class of representations of countable groups is a
natural extension of the class of representations with finite traces,
i.e., representations determined by finite characters? Recall that a
character of a group $G$ is a positive definite central 
(i.e., such that $\chi(gh)=\chi(hg)$ for all $g,h \in G$) function on $G$
whose value at the group identity is equal to one. According to the
GNS construction, 
an indecomposable character (which cannot be written as a nontrivial convex
combination of other characters) determines  a factor representation
(f.r.)\ of finite type
${\rm I}_n$, $n<\infty$, or ${\rm II}_1$.
Conversely, a f.r.\ of finite type uniquely determines a
character. For many (though not all) countable groups, 
the set of all f.r.\ of type~${\rm II}_1$ 
has the structure of a standard Borel space. At the same
time, for countable groups 
that have no
Abelian subgroup of finite index, the set of classes of pairwise
nonequivalent irreducible representations is not tame. In particular,
it has
no standard Borel structure. Hence for these groups, there is no
reasonable classification of irreducible representations. 
Moreover, there exist representations that have several different
decompositions into irreducible components. Such groups are usually
called wild. The infinite symmetric group 
$\mathfrak{S}_{\mathbb{N}}$, which consists of all finite permutations
of the set of integers ${\mathbb{N}}$, is a typical example of a wild
group. However, its f.r.\
of type ~${\rm II}_1$ are completely classified in
\cite{Thoma,VK0,Ok} and
are used for constructing harmonic analysis. 
A similar situation holds
for the infinite-dimensional unitary group~$\mathbb{U}(\infty)$. At the same time, important classes of countable
groups, in particular, the infinite general linear group over a finite
field $GL(\infty,\mathbb{F}_q)$, have too few finite characters for constructing a nontrivial
harmonic analysis. For instance, the characters
of $GL(\infty,\mathbb{F}_q)$ do not separate
points. These
considerations lead to the need for an extension of the class of
representations of finite type, i.e., for a generalization of the notion
of a character.

But the class of all representations of  type~${\rm II}_{\infty}$ is
too wide: e.g., for the infinite symmetric group
$\mathfrak{S}_{\mathbb{N}}$, it is as wild (and in the same
sense) as the set of classes of irreducible representations.
Therefore, a putative extension must include some reasonable restrictions
on representations of  type~${\rm II}_{\infty}$. An attempt to construct
such an extension was
made in \cite{VK0}, where a
class of representations 
 determined by semifinite traces on the group algebra was introduced. But
this class turned out to be too narrow, it does not even cover all
representations of 
$\mathfrak{S}_\mathbb{N}$  of type~${\rm II}_{\infty}$
associated with admissible representations of
$\mathfrak{S}_\mathbb{N}\times\mathfrak{S}_\mathbb{N}$.

{\it In this paper, we introduce the class of stable representations, which
is a natural extension of the class of representations of finite type.} 
It turns out that
stable representations of~$\mathfrak{S}_{\mathbb{N}}$ are of type~${\rm II}_{\infty}$. 
We also give a complete classification of stable f.r.\
up to quasi-equivalence. At the same time, we obtain
an answer to the question posed by the first author \cite{VKOk} in
connection with G.~I.~Olshanski's \cite{Ol2} theory of admissible
representations of the group
$\mathfrak{S}_\mathbb{N}\times\mathfrak{S}_\mathbb{N}$, that of
identifying the components of an admissible representation.
Namely, we prove that
the set
of stable f.r.\ coincides with the class of
representations that can be obtained as the restrictions of admissible irreducible
representations of $\mathfrak{S}_\mathbb{N}\times\mathfrak{S}_\mathbb{N}$
to the left and right components
($\mathfrak{S}_{\mathbb{N}}\times e$  
and $e\times\mathfrak{S}_{\mathbb{N}}$, respectively).

\section{Basic notions}

We introduce a new notion of stability for a positive
definite function (p.d.f.)\ on a group and the corresponding representation.
Typical groups for which this notion is meaningful and amenable to study are
inductive limits of compact groups.

\subsection{Topology on groups of automorphisms of a group}

Let $G$ be a countable group and
$\mathbf{C}^\star(G)$ be its group
$\mathbf{C}^\star$-algebra. We identify $G$ with
its natural image in $\mathbf{C}^\star(G)$, and
positive functionals (p.f.)\ on
$\mathbf{C}^\star(G)$ with the corresponding p.d.f.\ on~$G$. A
p.f.\ is called a state if it is equal to one at
the group identity. Let
${\rm Aut}\,G$ be the group of all automorphisms of~$G$. An element
$g\in G$ determines an inner automorphism
${\rm Ad}\,g\in {\rm Aut}\,G$: ${\rm
Ad}\,g(x)=gxg^{-1}$, $x\in G$. The group ${\rm Int}\,G$ of inner automorphisms is a
normal subgroup  in ${\rm Aut}\,G$.

We endow ${\rm Aut}\,G$ with the strong topology, in which a base of
neighborhoods of the identity automorphism consists of the sets
\begin{eqnarray}\label{strong_discrete_top}
\mathcal{U}_g=\left\{\theta \in {\rm Aut}\,G:  \theta(g)=g\right\}, \;\;\; g\in G.
\end{eqnarray}
Now we introduce the strong topology on  ${\rm Aut}\,G$  for an arbitrary locally compact
(l.c.)\ group. Let
${\rm Aut}\,\mathbf{C}^\star(G)$ be the group of all
automorphisms of the algebra
$\mathbf{C}^\star(G)$. We identify
${\rm Aut}\,G$ and ${\rm Int}\,G$ with the corresponding subgroups in
${\rm Aut}\,\mathbf{C}^\star(G)$. A base of the {\it strong topology} on
 ${\rm Aut}\,\mathbf{C}^\star(G)$ is determined by the neighborhoods
 $\mathcal{U}_{a,\epsilon }$ of
 the identity automorphism, where
\begin{eqnarray}\label{strong_aut_topology}
\mathcal{U}_{a,\epsilon }= \left\{\theta\in{\rm
Aut}\,\mathbf{C}^\star(G):\left\|\theta(a)-a \right\|<\epsilon
\right\},\;\;\;a\in \mathbf{C}^\star(G),\; \epsilon >0.
\end{eqnarray}
Then ${\rm Aut}\,G$ is a closed subgroup in
${\rm Aut}\,\mathbf{C}^\star(G)$. Denote by $\overline{\rm Int}\, G$
the completion of ${\rm Int}\,G$.

Let $G=\lim G_i$ be the inductive limit of a sequence of 
l.c.\ groups
$\left\{ G_i\right\}_{i\in \mathbb{N}}$, where
$G_i$ is a closed subgroup in
$G_{i+1}$ for all $i$. A base of the strong topology on
${\rm Aut}\,G$ is determined by the neighborhoods
$\mathcal{U}_{n,a,\epsilon }$, $a\in
\mathbf{C}^\star\left(G_n \right)$, where
\begin{eqnarray}\label{strong_aut_topology_ind}
\mathcal{U}_{n,a,\epsilon }= \left\{\theta\in {\rm Int}\,G:\theta\left(G_n
\right)=G_n \text{ and }\left\|
\theta(a)-a\right\|_{\mathbf{C}^\star\left(G_n
\right)}<\epsilon\right\}.
\end{eqnarray}

\medskip
{\bf 1.}\label{exemple_one} Our basic example is as follows.
Let $\mathbb{N}$ be the set of positive integers. A bijection
$s: \mathbb{N}\to \mathbb{N}$ is called {\it finite}
if the set $\left\{ i\in\mathbb{N}\big|\,s(i) \neq
i\right\}$ is finite. We define
$\mathfrak{S}_\mathbb{N}$ as the group of all finite bijections
$\mathbb{N}\to\mathbb{N}$ and set
$\mathfrak{S}_n=\left\{s\in\mathfrak{S}_\mathbb{N}|\; s(i)=i\right.$
for all $\left.
i>n\right\}$. Denote by $\mathfrak{S}_{\mathbb{N}\setminus n}$
the subgroup in  $\mathfrak{S}_\mathbb{N}$ consisting of the
elements leaving the
numbers $1,2,\ldots,n$ fixed  ($n<\infty$).

If $\overline{\mathfrak{S}}_\mathbb{N}$ is the group of all
bijections of $\mathbb{N}$, then
 $\mathfrak{S}_\mathbb{N}\subset\overline{\mathfrak{S}}_\mathbb{N}$,
 and for every $s\in\overline{\mathfrak{S}}_\mathbb{N}$, the map
$\mathfrak{S}_\mathbb{N} \ni x\mapsto sxs^{-1}\in\mathfrak{S}_\mathbb{N}$
is an automorphism ${\rm Ad}\, s$
of $\mathfrak{S}_\mathbb{N}$, which can be
naturally extended to an automorphism of the group
$\mathbf{C}^\star$-algebra
$\mathbf{C}^\star\left(\mathfrak{S}_\mathbb{N} \right)$ of
$\mathfrak{S}_\mathbb{N}$. One can easily check that ${\rm
Ad}\,\overline{\mathfrak{S}}_\mathbb{N}$ coincides with ${\rm
Aut}\,\mathfrak{S}_\mathbb{N}$  and is the closure of ${\rm
Int}\,\mathfrak{S}_\mathbb{N}$ in each of the topologies
(\ref{strong_discrete_top}), (\ref{strong_aut_topology}), and
(\ref{strong_aut_topology_ind}), which are equivalent. 
In particular, 
$\mathfrak{S}_\mathbb{N}$ is a dense normal subgroup in
$\overline{\mathfrak{S}}_\mathbb{N}$.

\medskip
{\bf 2.} Consider the infinite-dimensional unitary group
$U(\infty)$, which is the inductive limit $U(\infty)=\lim\limits_\rightarrow\coprod\limits_n U(n)$
of the finite-dimensional unitary groups
with respect to the natural
embeddings. If ${\rm Ad}\, g\in\mathcal{U}_{n,a,\epsilon }$, $g\in
U(\infty)$, then
$g=g(n)\cdot g_\infty(n)$,
where
$g(n)\in U(n)$ and $g_\infty(n)\cdot u=u\cdot g_\infty(n)$ for all $u\in
U(n)$. It is not difficult to show that
$\min\limits_{z\in\mathbb{C}:|z|=1}\left\| g(n)-zI_n\right\|\to 0$ as $\epsilon \to 0$
and ${\rm Int}\, U(\infty)\neq \overline{\rm
Int}\,U(\infty)$.

\subsection{Characters and representations of  $\overline{\rm
Int}\,G$.}\label{rep_Int}

Let $G$ be a countable group, and let 
$\Pi$ be the biregular representation\footnote{Given by the formula 
$\left( \Pi(g,h)\eta\right)(x)=\eta\left( g^{-1}xh\right) $, $\eta\in l^2(G)$.}
of the group $G\times G$. Then it follows from~(\ref{strong_aut_topology})
that the map
${\rm Int}\,G\ni{\rm Ad}\,g\stackrel{\mathcal{A}}{\mapsto}\Pi((g,g))\in \mathbb{U}(l^2)
$
can be extended by continuity to
$\overline{\rm Int}\, G$. 
This continuity is preserved if we replace the biregular representation of
 $G$ by a representation corresponding to a character.

On the other hand, in the next section we will describe a construction
of a unique, up to unitary
equivalence, representation of 
$G\times G$, which plays the role of the biregular representation, for an arbitrary
p.d.f.\ on~$G$. But
the corresponding representation of
${\rm Int}\,G$ is,  in general, no longer  continuous. Hence it
cannot be extended by continuity to
$\overline{\rm Int}\, G$. However, there is a class of
states on $G$, containing finite characters, for which this
continuity persists. This is exactly the class of stable
states.

\subsection{The canonical construction of representations of
$G\times G$ and ${\rm Int}\,G$}\label{canonoc}
Let $\pi$ be the representation of a group $G$ corresponding to
a state $\varphi $ and $M$ be the $w^*$-algebra generated by the
operators from $\pi(G)$. Let $M^\prime$ be the commutant of~$M$,  $M_*^+$ 
be the set of weakly continuous positive functionals on $M$, 
$\widetilde{\varphi}\in M_*^+$ be a faithful state, and
${\theta}$ be the representation of $M$ corresponding to
$\widetilde{\varphi}$ in a Hilbert space
$\widetilde{H}$ with a bicyclic vector $\widetilde{\xi }$. Since
$\widetilde{\varphi}$ is faithful, the map
$M\stackrel{{\theta}}{\mapsto}\theta(M)=\widetilde{M}$ is an isomorphism
(see \cite{TAKES}).

Denote by $\widetilde{\mathcal{J}}$ the antilinear isometry from the
polar decomposition of the closure of the map
$\widetilde{M}\widetilde{\xi }\ni
m\widetilde{\xi }\stackrel{\widetilde{S}}{\mapsto}m^*\widetilde{\xi }$ $\in
\widetilde{M}\widetilde{\xi }$, $m\in\widetilde{M}$.
Set
$
\tilde{\pi}=\tilde{\theta}\circ\pi
$.
Since
$\widetilde{\mathcal{J}}\widetilde{M}\widetilde{\mathcal{J}}^{-1}=\widetilde{M}^\prime$,
the map $G\times G\ni (g,h)$
$\mapsto
\widetilde{\Pi}((g,h))=\widetilde{\pi}(g)\widetilde{\mathcal{J}}
\widetilde{\pi}(h)\widetilde{\mathcal{J}}^{-1}
$ is a representation.

Since $\widetilde{\mathcal{J}}m\widetilde{\mathcal{J}}^{-1}=m^*$ for all $m$
from the center of the algebra $\widetilde{M}$ (see \cite{TAKES}), a
representation $\mathcal{A}_\pi$ of the group ${\rm Int}\,G$ is well
defined by analogy with Section~\ref{rep_Int}:
$
{\rm Int}\,G\ni{\rm
Ad}\,g\stackrel{\mathcal{A}_\pi}{\mapsto}\widetilde{\Pi}((g,g))
$.

\subsection{Stable states and representations; admissible
representations}\label{main_gen_stab}
Let $G$ be a l.c.\ group, $\mathfrak{A}=\mathbf{C}^\star(G)$, and 
$\mathfrak{A}^\ast$ be the dual space to $\mathfrak{A}$. Given a
state $\varphi\in\mathfrak{A}^\ast$, denote by $o_\varphi $ the map
${\rm Int}\,G\ni\theta\stackrel{o_\varphi}{\mapsto}\varphi \circ \theta\in\mathfrak{A}^\ast$.

\begin{Def}
A state $\varphi$ is called {\it stable} if the map
$o_\varphi$ is continuous in the strong topology on ${\rm Int}\, G$ (see
(\ref{strong_aut_topology})) and in the norm topology of the dual space on
    $\mathfrak{A}^\ast$.
 \end{Def}

\begin{Rem}
For every $\psi\in\mathfrak{A}^\ast$, the map $o_\psi$ is continuous
in the weak topology of the dual space on $\mathfrak{A}^\ast$.
  \end{Rem}

The representation of $G$ corresponding to a stable state will
also be called {\it stable}.

To define a stable p.d.f.\ on an inductive limit
$G=\lim_{n}G_n$ of l.c.\ groups, consider the cone
$C^+_G$ of p.d.f.\ on $G$ with the
topology defined by the metric
$\rho(\psi,\varphi ) =\sup\limits_n  \left\| \psi - \varphi\right\|_n$,
$\varphi ,\psi\in C^+_G$,
where $\|\cdot\|_n$ stands for the norm of the dual space
$\mathbf{C}^\star( G_n)^\ast$.

\begin{Def}
A p.d.f.\ $\varphi$ on $G$ is called stable if the map
${\rm Int}\,G\ni\theta\stackrel{o_\varphi}{\mapsto}\varphi \circ \theta\in C^+_G$
is continuous in the topology on
${\rm Int}\, G$ defined according to~(\ref{strong_aut_topology_ind}).
\end{Def}

\begin{Rem} If $G$ is a free group, then the
topology~(\ref{strong_discrete_top}) on ${\rm Aut}\, G$
is discrete. Hence all p.d.f.\ on $G$ are stable.
\end{Rem}

Now we introduce the 
notion of an admissible representation of 
$G\times G$, which extends the corresponding notion from
\cite{Ol2} to the case of an arbitrary group.

 \begin{Def}
Let  $G$ be a l.c.\ group or an inductive limit of l.c.\
groups. A unitary representation
$\Pi$ of the group $G\times G$ in a Hilbert space $H$ is called
admissible if the map
 $
 {\rm Int}\,G\ni{\rm Ad}\,g$ $\mapsto \Pi((g,g))\in \mathbb{U}(H)
 $
is continuous in the strong topology on
${\rm Int}\,G$ (see~(\ref{strong_aut_topology}), (\ref{strong_aut_topology_ind}))
and in the strong operator topology on
$\mathbb{U}(H)$.
 \end{Def}

\begin{Th}
If the representation $\widetilde{\Pi}$ from Section~{\rm\ref{canonoc}}    
corresponds to a stable representation $\pi$ of $G$, then 
$\widetilde{\Pi}$ is an admissible representation of the group $G\times G$.   
The restriction of $\widetilde{\Pi}$ to $\mathfrak{S}_\mathbb{N}\times
e$ is quasi-equivalent to $\pi$, the center
$\mathcal{C}$ of the algebra 
$\widetilde{\Pi}\left(\mathfrak{S}_\mathbb{N}\times\mathfrak{S}_\mathbb{N}\right)
^{\prime\prime}$  coincides with the center of $\widetilde{M}$,   
and the components of the decomposition of $\widetilde{\Pi}$ with
respect to $\mathcal{C}$ are irreducible representations. 
\end{Th}

\section{Stable representations of the group
$\mathfrak{S}_\mathbb{N}$ and their classification}
In this section, we formulate results on the classification of stable
representations of 
$\mathfrak{S}_\mathbb{N}$ up to quasi-equivalence.

With a stable f.r.\
$\pi$ of the group $\mathfrak{S}_\mathbb{N}$ we associate an invariant
$\chi^{\pi\,a}$, called an {\it asymptotic character}. Let
$M=$ $\pi\left(\mathfrak{S}_\mathbb{N}\right)^{\prime\prime}$.

\begin{Prop}\label{as_character}
Let $\pi$ be a stable f.r.\ of
$\mathfrak{S}_\mathbb{N}$. For every $g\in\mathfrak{S}_\mathbb{N}$ 
there exists a sequence
$\left\{ \sigma_n^g
\right\}_{n\in\mathbb{N}}\subset\mathfrak{S}_\mathbb{N}$ such that
$\sigma_n^gg\left( \sigma_n^g\right)^{-1}\in \mathfrak{S}_{\mathbb{N}\setminus n}$
and $\sigma_{n+1}^g \left( \sigma_n^g\right)^{-1}
\in\mathfrak{S}_{\mathbb{N}\setminus n}$. If $\psi\in
M_*^+$, then the limit 
$\lim\limits_{n\to\infty}\psi\left( \pi\left(\sigma_n^gg\left(
\sigma_n^g\right)^
{-1} \right)\right)\cdot\psi(I)^{-1}$ $=\chi^{\pi\,a}(g)$
exists and does not depend on
$\psi$ and  $\left\{ \sigma_n^g \right\}_{n\in\mathbb{N}}$. The
function $\chi^{\pi\,a}$ is an indecomposable character of
$\mathfrak{S}_\mathbb{N}$ (see~{\rm\cite{Thoma}}). If $\tilde{\pi}$
is quasi-equivalent to $\pi$,  
then  $\chi^{\pi\,a}$ $=\chi^{\tilde{\pi}\,a}$.
\end{Prop}

\begin{Th}\label{stab_centr}
Fix a representation $\pi$ of
$\mathfrak{S}_\mathbb{N}$, and let\linebreak $M_{*}^+(n)= \left\{ \phi\in M_*^+: \phi \left(\pi(s) x\right)=
\phi(x\pi(s))\mbox{ for all } x\in M, s\in\mathfrak{S}_{\mathbb{N}\setminus
n}\right\}
$. The following conditions are equivalent:
\begin{itemize}
 \item[\rm i)] the representation $\pi$ is stable;
 \item[\rm ii)] the union $\bigcup\limits_{n}M_*^+(n)$ is dense in
 the norm topology of the dual space on~$M_*^+$.
 \end{itemize}
\end{Th}

We define  the {\it central depth} ${\rm cd}(\pi)$ of a representation
$\pi$ as ${\min}\left\{n: M_*^+(n)\neq0 \right\}$.
Clearly, ${\rm cd}(\pi)$ is a quasi-equivalence
invariant.

\begin{Th}\label{strong_ind}
Let $\pi$ be a stable f.r.,
$n={\rm cd}(\pi)$, and
$\psi\in M_*^+(n)$. If $E$ is the support\footnote{That is, $E$ 
is the smallest projection from 
$\pi\left(\mathfrak{S}_\mathbb{N}\right)^{\prime\prime}$ such that $\psi({\rm
I}-E)=0$.} of $\psi$, then $E\in \pi\left(\mathfrak{S}_{n}\mathfrak{S}_{\mathbb{N}\setminus{n}}
\right)^{\prime}$ and $E\pi(s)E=0$ for all
$s\notin\mathfrak{S}_{n}\mathfrak{S}_{\mathbb{N}\setminus{n}}$.
In particular, the representation
$\pi_E$ of the group $\mathfrak{S}_{n}\mathfrak{S}_{\mathbb{N}\setminus n}$
determined by the operators
$E\pi(s)E$ has the form
 $T_\lambda\otimes \pi_{\alpha \beta }$, where $\lambda\vdash n$ 
is a Young diagram,
$T_\lambda$ is the corresponding irreducible representation of
$\mathfrak{S}_n$, and $\pi_{\alpha \beta }$ is the representation of
$\mathfrak{S}_{\mathbb{N}\setminus n}$ corresponding to a Thoma
character $\chi_{\alpha\,\beta}$. Hence $\pi$ is 
quasi-equivalent to the representation
${\rm Ind}_{\mathfrak{S}_{n}
\mathfrak{S}_{\mathbb{N}\setminus n}}^{\mathfrak{S}_\mathbb{N}}T_\lambda \otimes
\pi_{\alpha \beta }$.
\end{Th}

The next result describes the dependence of ${\rm cd}(\pi)$
on the Thoma parameters
${\alpha, \beta }$ of the asymptotic character $\chi^{\pi\,a}$.

\begin{Th}\label{ind}
Let $n\geq 1$, and let $\lambda$, $T_\lambda$, $\pi_{\alpha \beta }$ be
as in Theorem~{\rm\ref{strong_ind}}. Then $\Pi_{\alpha \beta }^\lambda$ is
a f.r.\ of
type~${\rm II}_\infty$ if $\sum\alpha_i +\sum\beta_i
=1$,  and a f.r.\ of
type~${\rm II}_1$ if
$\sum\alpha_i
+\sum\beta_i$ $<1$.
\end{Th}

\begin{Th}\label{nonquasi}
Let 
$\sum\alpha_i +\sum\beta_i
=1$. The representations $\Pi_{\alpha \beta }^\lambda$ and
$\Pi_{\tilde{\alpha}\tilde{\beta} }^{\tilde{\lambda}}$ are
quasi-equivalent if and only if
$\alpha_i=\tilde{\alpha}_i $,  $\beta _i=\tilde{\beta}_i$, and $\lambda=\tilde{\lambda}$.
\end{Th}

Now we describe a bijection between the set of quasi-equivalence
classes of stable f.r.\ and the set ``partially central'' states on 
$\mathfrak{S}_\mathbb{N}$.

\begin{Th}\label{unique}
Let $\pi$ be a stable f.r.\ of
$\mathfrak{S}_\mathbb{N}$ and $n={\rm cd}(\pi)$.
If a state $f$ on
$\mathfrak{S}_\mathbb{N}$ determines a representation
quasi-equivalent to $\pi$ and satisfies the  condition
 $f\left(tst^{-1} \right)=f\left(s \right)$ for
all $s\in\mathfrak{S}_\mathbb{N}$ and
$t\in\mathfrak{S}_n\mathfrak{S}_{\mathbb{N}\setminus n}$, then
\begin{eqnarray}
f(s)
   =\left\{
\begin{array}{rl}
\chi_\lambda\left(s_1 \right)
 \chi_{\alpha \,\beta }\left(s_2 \right)& \text{ if }
 s_1\in\mathfrak{S}_n,\, s_2\in\mathfrak{S}_{\mathbb{N}\setminus n},
 \text{ and } s=s_1s_2,\\
0 & \text{ if } s\notin
 \mathfrak{S}_n\cdot \mathfrak{S}_{\mathbb{N}\setminus n},
 \end{array}\right.
\end{eqnarray}
where $\lambda\vdash n$ is a diagram, 
$\chi_\lambda$ is the normalized character of the corresponding
irreducible representation of $\mathfrak{S}_n$, 
and  $\alpha,\beta$
are the Thoma parameters of the asymptotic character $\chi^{\pi\,a}$
 (see Proposition~{\rm\ref{as_character}}).
\end{Th}

\medskip
Translated by N.~V.~Tsilevich.

{}

\end{document}